\documentclass[12pt,leqno]{article}
\usepackage{color}
\definecolor{purple}{rgb}{0.65, 0, 1}
\definecolor{orange}{rgb}{1,.5,0}
\definecolor{brown}{rgb}{.9,.73,.26}

\usepackage{amsfonts}
\usepackage{amsmath}
\usepackage{epsfig}
\usepackage{graphicx}
\input epsf

\def\div{\hbox{div}}
\def\R{\hbox{\bf R}}

\def\div{\hbox{div}}

\def\N{\hbox{\bf N}}

\def\e{\epsilon}

\def\d{\delta}

\newcommand{\ba}{\begin{eqnarray}}
\newcommand{\ea}{\end{eqnarray}}

\newtheorem{theo}{\bf Theorem}[section]
\newtheorem{lem}[theo]{\bf Lemma}
\newtheorem{pro}[theo]{\bf Proposition}
\newtheorem{cor}[theo]{\bf Corollary}
\newtheorem{defi}[theo]{\bf Definition}
\newtheorem{rem}[theo]{\bf Remark}

\numberwithin{equation}{section}

\renewcommand{\N}{{\mathbb N}}
\renewcommand{\R}{{\mathbb R}}
\renewcommand{\P}{{\mathcal P}}

\renewcommand{\e}{\varepsilon}

\newenvironment{Proofc}[1]{\smallskip\par\noindent\textsc{#1}\quad}%
  {\hfill$\Box$\bigskip\par}

\newcommand{\loc}{{\textrm{loc}}}

\newcommand{\<}{\langle}
\renewcommand{\>}{\rangle}

\newcommand{\sgn}{\operatorname{sgn}}

\renewcommand{\P}{\mathcal{P}}
\newcommand{\sphere}{\mathbb{S}}

\newcommand{\dd}{\partial}
\renewcommand{\d}{\,\textup{d}}

\newcommand{\supp}{\operatorname{supp}}

\newcommand{\lap}{\Delta}

\renewcommand{\div}{\operatorname{div}}

\begin{document}

\title{\bf Pointwise estimates for the heat equation. 
Application to the free boundary of the obstacle problem
with Dini coefficients}
 \author{E. Lindgren\footnote{{erik.lindgren@math.ntnu.no, Dept. of Mathematical
Sciences, NTNU, 7491 Trondheim, Norway}} and R. Monneau\footnote{Universit\'e Paris-Est, CERMICS, Ecole des Ponts ParisTech,
   6 et 8 avenue Blaise Pascal, Cit\'e Descartes, Champs-sur-Marne,
   77455 Marne-la-Vall\'ee Cedex 2}}
\maketitle


\vspace{20pt}


 \centerline{\small{\bf{Abstract}}}
\noindent {\small{We study the pointwise regularity of
    solutions to parabolic equations. As a first result, we prove that if
    the modulus of mean oscillation of $\Delta u -u_t$ at the origin is Dini (in $L^p$
    average), then the origin is a Lebesgue point of continuity
    (still in $L^p$ average) for $D^2 u$ and $\dd_t u$. 
    We extend this pointwise regularity
    result to the parabolic obstacle problem with Dini
    right hand side. In particular, we prove that the
    solution to the obstacle problem has, at regular points of the free boundary, a Taylor expansion up to order
    two in space and one in time (in the $L^p$ average). 
    Moreover, we get a quantitative estimate of
    the error in this Taylor expansion.
    Our method is based on decay estimates obtained by
    contradiction, using blow-up arguments and Liouville type theorems.
    As a by-product of our approach, we deduce that the regular points of the free boundary are locally contained in a $C^1$ hypersurface for the parabolic distance $\sqrt{x^2 +|t|}$.
}}\hfill\break

 \noindent{\small{\bf{AMS Classification:}}} {\small{35R35.}}\hfill\break
 \noindent{\small{\bf{Keywords:}}} {\small{Obstacle problem, Heat
     equation, Dini condition, free boundary, pointwise regularity.
}}\hfill\break


\section{Introduction}

\subsection{The heat equation}

In this paper, we are interested in the pointwise regularity of solutions to parabolic problems. We first consider the solutions to the following heat equation
\begin{equation}\label{eq:100}
\left\{\begin{array}{l}
\Delta u -u_t = f \quad \mbox{in}\quad Q_1^-,\\
\\
f\in  L^p\left(Q_1^-\right) \quad \mbox{and}\quad f(0)=0,
\end{array}\right.
\end{equation}
where we set the past cylinder $Q_r^-=B_r\times (-r^2,0]$ with $B_r=B_r(0)$ the open ball in $\R^n$, of radius $r$ centered at the origin $0$. Here $p\in (1,+\infty)$ and we assume that $0=(0,0)$ is a Lebesgue point of $f$, in order to define $f(0)$ if necessary.\\

It is well-known that if $f$ is H\"older continuous in the cylinder $Q_1^-$, then so are the spatial second derivatives of $u$ and the time first derivative of $u$ (see for instance \cite{Lie96}).
Let us introduce the following parabolic modulus of continuity of $f$ on the cylinder $Q_1^-$
$$
\overline{\sigma}(r)=\sup_{\stackrel{|x-y|^2 + |t-s|\le r^2}{ (x,t),(y,s)\in Q_1^-}} |f(x,t)-f(y,s)|.
$$

\begin{defi}{\bf (Dini function)}
\\
A function $\overline{\sigma}$ is said to be Dini if
\begin{equation*}
\int_0^1\frac{\overline{\sigma}(r)}{r}\ dr < +\infty.
\end{equation*}
\end{defi}
It is well-known (see  \cite{W06}) that if $\overline{\sigma}$ is Dini, then 
the second derivatives of $u$ are continuous in the cylinder $Q_{1/2}^-$ 
with a modulus of continuity proportional to 
$$r \sup_{Q_1^-}|u| + \int_0^{r}\frac{\overline{\sigma}(s)}{s}ds +  r \int_r^{1}\frac{\overline{\sigma}(s)}{s^2}ds.$$
Notice that the modulus of continuity of $u_t$ then follows from equation (\ref{eq:100}) itself.\\

Up to our knowledge, such results are usually obtained assuming a modulus of continuity in an open set.
Here we change the point of view, and only consider {\it pointwise modulus of mean oscillation},
like for instance \cite{ZC02}.
For any $p\in (1,+\infty)$, we define a kind of modulus of mean oscillation (in $L^p$ average) of the
function $f$ {\it at the
origin} as
\begin{equation}\label{eq:sigma-tilde}
\tilde{\omega}(r) =\tilde{\omega}(f,r)=\ \inf_{c\in\R} \ \left(\frac{1}{|Q_r^-|}\int_{Q_r^-}
    |f(x,t)-c|^{p}\right)^{\frac1{p}}.
\end{equation}
Furthermore, we denote by $\tilde{\mathcal P}_2$ the set of polynomials of degree
less than or equal to two in space and of degree less than or equal to one in time. Let
\begin{equation}\label{eq::e11}
\tilde{N}(u,\rho)=\inf_{P\in \tilde{\P}_2}\left(\frac{1}{r^{n+2+2p}}\int_{Q_r^-}
  |u-P|^p\right)^{\frac1p}.
\end{equation}

\begin{theo}\label{th:1-bis}{\bf (Pointwise parabolic BMO estimates for the heat equation)}\\
Let $p\in (1,+\infty)$. Then there exist $\alpha\in (0,1]$
and constants $r_*\in (0,1]$, $C>0$, such that the following holds.
If $u\in L^p(Q_1^-)$ satisfies (\ref{eq:100}) with
the associated $\tilde{\omega}$ defined in
(\ref{eq:sigma-tilde}), then we have:\\
\noindent {\bf i) Pointwise BMO estimate}
\begin{equation}\label{eq::e25}
\sup_{r\in (0,1]}\tilde{N}(u,r)  \le C\left\{ \left(\int_{Q_1^-}
    |u|^p\right)^{\frac{1}{p}}+\left(\int_{Q_1^-}
    |f|^p\right)^{\frac{1}{p}}+ \sup_{r\in (0,1]}\tilde{\omega}(r)\right\}   
\end{equation}
\noindent {\bf ii) Pointwise VMO estimate}\\
$$\displaystyle{\left(\tilde{\omega}(r)\longrightarrow 0  \quad \mbox{as}\quad r\to 0^+\right)
\quad \Longrightarrow \quad \left(\tilde{N}(u,r) \longrightarrow 0 \quad \mbox{as}\quad r\to 0^+\right)}$$
\noindent {\bf iii) Pointwise control on the solution}\\ 
If $\tilde{\omega}$ is {\bf Dini}, then $\tilde{N}(u,\cdot)$ is
{\bf Dini}, and there exists a caloric polynomial $P_0$ (i.e., a solution of $(P_0)_t=\Delta P_0$) of degree less than or
equal to two in space and of degree less than or equal to one in time, such that for every $r\in (0,r_*]$ there holds
\begin{equation}\label{eq::e26}
\left(\frac{1}{|Q_r^-|}\int_{Q_r^-}
    \left|\frac{u(x,t)-P_0(x,t)}{r^2}\right|^{p} 
  \right)^{\frac1{p}} \le C\left(\tilde{M}_0 r^\alpha + \int_{0}^r\frac{\tilde{\omega}(s)}{s}
    \, ds + r^\alpha \int_r^1 \frac{\tilde{\omega}(s)}{s^{1+\alpha}}  \, ds\right)
\end{equation}
and
$$P_0(x,t)= a + b\cdot x +\frac12 {}^t x\cdot c\cdot x + mt,$$
with
$$|a|+|b|+|c|+|m|\le C \tilde{M}_0 \quad \mbox{and}\quad \tilde{M}_0= \int_0^{1}
 \frac{\tilde{\omega}(s)}{s}\ ds 
 +  \left(\int_{Q_1^-} |u|^p\right)^{\frac1p} + \left(\int_{Q_{1}^-} |f|^p\right)^{\frac1p}. $$
\end{theo}


\begin{rem}
Theorem  \ref{th:1-bis} iii)  implies in particular (using parabolic estimates) that we have
a Lebesgue point of continuity of the second derivatives  $D^2u$ and of $u_t$ (in the $L^p$ average) if $\Delta u-u_t$
has a Dini modulus of mean oscillation (in the $L^p$ average) at the same point.
\end{rem}

\begin{rem}
A straightforward consequence of Theorem    \ref{th:1-bis} is
that the second derivatives $D^2u$ and $u_t$ are H\"{o}lder continuous in an open set
$\Omega\subset \R^{n+1}$, if $\Delta u-u_t$ is H\"{o}lder continuous in $\Omega$ 
for the parabolic distance $\sqrt{|x|^2+|t|}$.
\end{rem}

\begin{rem}
Notice that our definition of $\tilde{\omega}(r)$ differs from the analogue given in \cite{Mon09},
not only because we consider here the parabolic problem instead of the elliptic one, 
but also because there is no supremum in this new definition. 
From that point of view, estimate (\ref{eq::e26})
is finer than the one given in \cite{Mon09}, and than the ones that can be found in the classical literature.
\end{rem}

We would like to emphasize that the result of Theorem \ref{th:1-bis} is completely {\it pointwise}, which
does not seem to be so usual in the literature.

\subsection{The model obstacle problem}

In the second part of this article we are in particular interested in the regularity of the free boundary
for solutions to the parabolic obstacle problem. The model problem is the following. Consider a
function $u$ satisfying
\begin{equation}\label{eq:obstacle}
\left\{
\begin{array}{l}
\left. \begin{array}{l}
\Delta u -u_t= f(x,t)\chi_{\{u>0\}}   \\
\\
\displaystyle{u\geq 0}
\end{array}  \right| \quad\mbox{in}\quad Q_1^-,\\
\\
u, f \ \in L^p(Q_1^-)\quad \mbox{and}\quad f(0)=f(0,0)=1,\\
\\
0\in \partial \left\{u>0\right\},
\end{array} \right.
\end{equation}
for $p\in ((n+2)/2,+\infty )$, where $Q_1^-$ is the past unit cylinder as before and $\chi_{\{u>0\}}$ is the characteristic function of the set
$\left\{u>0\right\}$, which is equal to $1$ if $u>0$ and $0$ if $u=0$. From classical parabolic estimates joint with Sobolev embeddings with our
assumption $p>(n+2)/2$, 
every solution $u$ is in particular continuous, which
allows us to consider the boundary of the open set
$\displaystyle{\left\{u>0\right\}}$. Here $\partial\left\{u>0\right\}$ is
called the free boundary.
Moreover, we assume that 
$(0,0)$ is a Lebesgue point for $f$ in order to define $f(0)$.\\

There is a vast literature on the above problem. In the special case when $f=1$ and in a slightly more general setting, it is proved in \cite{ASU00}, that the solution enjoys the optimal $C^{1,1}_x\cap C^{0,1}_t$ regularity. Moreover, in \cite{ASU03}, it is proved that the free boundary $\dd\{u>0\}$ is, close to the part of the fixed boundary where $u$ satisfies a homogeneous Dirichlet condition, the graph of a Lipschitz function. This was extended to a more general problem in \cite{CPS04}, where it is proved that the free boundary is at regular points a $C^1$ regular graph. Some partial regularity results are also proved in \cite{EL11}, under the assumptions that $f$ is H\"older continuous.

In the one dimensional setting, and under the assumption that $f$ is Dini continuous, there is a series of paper, \cite{BDM05}, \cite{BDM06} and \cite{Bla06}, where this problem is studied. There it is proved that the free boundary is $C^1$ regular at certain regular (see the next page) points, and also that the free boundary enjoys a certain structure at the other points, the so called singular points.

Let us introduce the following kind of {\it pointwise modulus of continuity} (in $L^p$ average) of the
function $f$ {\it at the origin}:
\begin{equation}\label{eq:sigma}
\sigma(r) =\sigma(f,r)=\sup_{0<\rho\le r} \omega(\rho)
\quad \mbox{with}\quad \omega(\rho)  =\omega(f,\rho)=\left(\frac{1}{|Q_\rho^-|}\int_{Q_\rho^-}
    |f(x,t)-f(0)|^{p}\right)^{\frac1{p}}  
\end{equation}
We have the following general regularity result.
\begin{pro}\label{pro:sing1}{\bf (Quadratic growth)}\\
Let $p\in ((n+2)/2,+\infty)$. Then there exists  a constant $C>0$ such that
if $u$ is a solution of (\ref{eq:obstacle}) with $\sigma$ bounded given by
(\ref{eq:sigma}), then
$$0\le  u(x,t)\le C_1 \left(|x|^2 +|t|\right)\textup{ in $Q_\frac12^-$},$$
where   $C_1= C \left(1+ \sigma(1)\right)$.
\end{pro}
In order to present our main result, we need to introduce the quantity
$$
M_{\textup{reg}}(u,r)=\sup_{\rho\in (0,r]}\left(\inf_{P\in \P_{\textup{reg}}}
\left(\frac{1}{\rho^{n+2+2p}}\int_{Q_\rho^-} |u-P|^p\right)^\frac{1}{p}\right),
$$
where 
$$
\P_{\textup{reg}}=\left\{P(x,t)=\frac{1}{2}(\max(0,x\cdot \nu))^2,
\nu\in  \sphere^{n-1}\right\}.
$$
Notice that as a consequence of Proposition \ref{pro:sing1}, $M_{\textup{reg}}(u,\rho)$ is bounded for $\rho\le 1/2$. Recall that if the  free boundary is smooth (or regular) around the origin, 
then it is known that the blow-up limit
of the solution (i.e., the limit of certain rescalings of the solution) at the origin is unique and is an element of the set $\P_{\textup{reg}}$.
Therefore we have in particular
\begin{equation}\label{eq:106}
\lim_{r\to 0^+} M_{\textup{reg}}(u,r)=0.
\end{equation}
More generally, we define the set of regular points as
$${\mathcal R}=\left\{\begin{array}{l}
(x_0,t_0)\in Q_1^-,\quad (x_0,t_0) \quad \mbox{is a Lebesgue point of $f$ with}\quad f(x_0,t_0)>0\\
\\
\quad \mbox{and}\quad \displaystyle \lim_{r\to 0^+} M_{\textup{reg}}\left(\frac{1}{f(x_0,t_0)}u(x_0+\cdot,t_0+\cdot),r\right)=0
\end{array}\right\}.$$
Our main result is the following:
\begin{theo}\label{th:main}{\bf (Modulus of continuity at a regular point of the free
  boundary)} \\
Let $p\in ((n+2)/2,\infty)$. There exist $\alpha\in (0,1]$ and constants
  $C>0$,$M_0,r_0\in (0,1)$ such that, given $u$ satisfying  (\ref{eq:obstacle}), we have the following property.
If the modulus of continuity $\sigma$ defined in (\ref{eq:sigma}) is assumed {\bf Dini}, and if
$$M_{\textup{reg}}(u,r_0)\leq M_0,$$
then there exists $P_0\in \P_\textup{reg}$ such that for every $r\in (0,r_0)$
$$
\left(\frac{1}{|Q_r^-|}\int_{Q_r^-}
    \left|\frac{u-P_0}{r^2}\right|^p\right)^\frac{1}{p}\leq
C\left(M_{\textup{reg}}(u,r_0)r^\alpha+\int_0^r\frac{\sigma(s)}{s}\d
  s+r^\alpha\int_r^1\frac{\sigma(s)}{s^{1+\alpha}}\d s\right).
$$
\end{theo}

\begin{rem}
 With the same methods, it would be possible to get a similar estimate for any $p\in (1,+\infty)$,
 but under the stronger assumption that the coefficient of the
 right hand side of the equation is bounded from above and from below,
 i.e., $0< \delta_0 \le f\le 1/\delta_0$.
\end{rem}
Remark that under the assumptions of Theorem \ref{th:main}, we recover in particular (\ref{eq:106}). As a corollary of Theorem \ref{th:main} and using a Weiss type monotonicity formula, 
we will show in a companion paper \cite{LM} the result below.

\begin{theo}{\bf (Regularity of the regular set of the free boundary, \cite{LM})}
\\
Consider a solution $u$  of (\ref{eq:obstacle}), and assume that $\overline{\sigma}$ defined in (\ref{eq:106}) is {\bf Dini} with $f\ge \delta_0>0$ on $Q_1^-$. Then for any point $(x_0,t_0)\in {\mathcal R}$, there exists a neighborhood $V$ of $(x_0,t_0)$ in $Q_1^-$, such that 
$V\cap \partial\left\{ u>0\right\}$ is locally a $C^1$ hypersurface with respect to the parabolic distance. More precisely,
up to a rotation of the spatial coordinates
$$V\cap \partial\left\{ u>0\right\} = \left\{(x,t) \quad \mbox{such that}\quad x_n=g(x',t) \quad \mbox{with}\quad  (x',t)\in V'\right\},$$
where $x'=(x_1,...,x_{n-1})$, the set $V'$ is an open set in $\R^n$, and $g: V'\to \R$ is a map satisfying
$$g(x'+h',t+k)=g(x,t)+ h'\cdot D_{x'}g(x',t) + o(\sqrt{(h')^2+|k|}),$$
with $D_{x'}g$ continuous on $V'$.
\end{theo}
In \cite{LM} we will also present a theory for the singular points of the free boundary, that is, for the complement of the regular part.


\subsection{Organization of the paper}
The organization of the paper is as follows. First, in Section \ref{sec:class},  we recall certain classical results concerning parabolic Sobolev spaces and parabolic equations. This is followed by Section \ref{sec:eqthm}, where we, by contradictory and blow-up type arguments, prove our main result for the heat equation, namely Theorem \ref{th:1-bis}.

In Section \ref{sec:growth} we turn our attention to the obstacle problem. We prove, using mainly standard techniques, quadratic growth estimates for the obstacle problem and in Section \ref{sec:nondeg}, we exploit a standard non-degeneracy result and obtain a related, somewhat more technical result, refered to as weak non-degeneracy. In the following section, namely Section \ref{sec:comp}, we provide a compactness result that we strongly use to prove the main theorem (Theorem \ref{th:main}) for the obstacle problem, which is proved, using contradictory and blow-up type arguments, in Section \ref{sec:decay}.
\subsection{Notation}
Throughout the whole paper we will use the notation below:\\
$$
\begin{array}{ll}
u_t=\dd_t u =\frac{\partial u}{\partial t}&\textup{- the time derivative}\\
\lap u=\sum_{i=1}^n \frac{\dd^2 u}{\dd {x_i}^2} &\textup{- the Laplace operator}\\
Hu:= \Delta u -u_t &\textup{- the heat operator}\\
Q_r(x_0,t_0)=B_r(x_0)\times (t_0-r^2,t_0+r^2) &\textup{- a parabolic cylinder}\\
Q_r^-(x_0,t_0)=B_r(x_0)\times (t_0-r^2,0]&\textup{- a half cylinder}\\
\partial_pQ_r^-(x_0,t_0) = \left((\partial B_r(x_0))\times [t_0-r^2,0)\right) \bigcup \left(B_r(x_0)\times \left\{0\right\}\right)&\textup{- the parabolic boundary}\\
Q_r = Q_r(0,0), \quad Q_r^-=Q_r^-(0,0), \quad \partial_p Q_r^- = \partial_p Q_r^-(0,0)&\textup{- simplified notation}\\
\omega(g,\rho)=\left(\frac{1}{|Q_\rho^-|}\int_{Q_\rho^-}
    |g(x,t)-g(0)|^{p}\right)^{\frac1{p}} & \textup{- the average oscillation over a cylinder}\\
 \sigma(g,r)=\sup_{0<\rho\le r}\omega(g,\rho) &\textup{- a special $L^p$-modulus}\\
\tilde \P_2& \textup{- polynomials of parabolic degree}\\
&\,\, \textup{   less than or equal to two}\\
\P_2& \textup{- $P\in \tilde \P_2$ such that $HP=0$}\\
 \P_{2,c}& \textup{- $P\in \tilde \P_2$ such that $HP=c$}
\end{array}
$$

\section{Classical results for parabolic equations}\label{sec:class}
Here we recall the following classical results that will be of constant use in the rest of the paper.
\begin{theo}\label{th:interior}{\bf (Parabolic interior $L^p$-estimates)} \\Let $p\in (1,\infty)$.
  If $u\in L^p(Q_r^-)$ and $Hu\in L^p(Q_r^-)$ then 
$$
||u||_{W^{2,1}_p(Q_{r/2}^-)}\leq
  C\left(||u||_{L^p(Q_r^-)}+||Hu||_{L^p(Q_r^-)}\right), 
$$
where
$$
W^{2,1}_p(Q_\rho^-)=\{v\in L^p(Q_\rho^-)\big|\quad  v,\nabla v,D^2 v,v_t\in
L^p(Q_\rho^-)\}, 
$$
endowed with the norm
$$
||u||_{W^{2,1}_p(Q_\rho^-)}=||u||_{L^p(Q_\rho^-)}+||\nabla
u||_{L^p(Q_\rho^-)}+||D^ 2u||_{L^p(Q_\rho^-)}+||u_t||_{L^p(Q_\rho^-)}.
$$
\end{theo}
The result above is a special case of Theorem 7.22 on page 175 in \cite{Lie96}.
\begin{theo}\label{th:sobolev} {\bf (Parabolic Sobolev embedding)}\\ Let $u\in W_p^{2,1}(Q_r^-)$ with $p\in ((n+2)/2,\infty)$. Then
$$
||u||_{C^{\alpha}(Q_r^-)}\leq C_*||u||_{W^{2,1}_p(Q_r^-)}, 
$$
with $\alpha=2-\frac{n+2}{p}$ and where $C^\alpha(Q_r^-)$ refers to the parabolic H\"older space.
\end{theo}
This result is contained in Lemma 3.3 on page 80 in \cite{LSU67}.

\begin{theo}\label{th:para} {\bf (Classical $L^p$ parabolic estimate)}\\ 
Let $u\in W_p^{2,1}(Q_r^-)$ for $p\in (1,+\infty)$ a solution of
$$\left\{\begin{array}{ll}
Hu = f &\quad \mbox{on}\quad Q_r^-,\\
u=0  &\quad \mbox{on}\quad \partial_p Q_r^-,
\end{array}\right.$$
where $\partial_p Q_r^-$ 
denotes the parabolic boundary of $Q_r^-$, and  $f\in L^p(Q_r^-)$. Then there exists a constant $C>0$ 
(depending only on $p$, the dimension $n$ and $r>0$) such that
$$||u||_{W^{2,1}_p(Q_r^-)}\le C ||f||_{L^p(Q_r^-)}.$$
\end{theo}
This result can be found in Proposition 7.18 on page 173 in \cite{Lie96}.

\section{Proof of Theorem \ref{th:1-bis}}\label{sec:eqthm}

In order to give the proof of Theorem \ref{th:1-bis}, we show a basic decay estimate in a first subsection and some routine results in a second subsection. The proof of Theorem \ref{th:1-bis} is done in the third subsection.

\subsection{A basic decay estimate}\label{sec:eqthm1}

Given a function $f$, we consider a (unique) constant $c_r$ such that
$$\tilde{\omega}(f,r)=\inf_{c\in \R} \left(\frac{1}{|Q_r^-|}\int_{Q_r^-} |f(x,t)-c|^p\right)^{\frac1p} 
= \left(\frac{1}{|Q_r^-|}\int_{Q_r^-} |f(x,t)-c_r|^p\right)^{\frac1p}.$$
We define the particular set of caloric polynomials:
$$ \P_2=\left\{P \mbox{ a caloric polynomial}
\quad \left|\begin{array}{l}
\mbox{of degree less than or equal to $2$ in space}\\
\\
\mbox{of degree less than or equal to $1$ in time}
\end{array}\right.\right\}$$
Considering a  particular polynomial $P_*\in \tilde{\P}_2$ which satisfies $\Delta P_* -(P_*)_t =1$
(for instance $\displaystyle P_*(x,t)=\frac{x^2}{2n}$), we define
$$\P_{2,c}= c P_* + \P_2$$
and for a function $u$ solving (\ref{eq:100}) we let
\begin{equation}\label{eq::e15}
\hat{N}(u,r)=\inf_{P\in \P_{2,c_r}}\left(\frac{1}{r^{n+2+2p}}\int_{Q_r^-} |u-P|^p\right)^{\frac1p}.
\end{equation}
For $\tilde{\omega}(f,r)$ and $\hat{N}(u,r)$ respectively defined in (\ref{eq:sigma-tilde}) and (\ref{eq::e15}),
we now define for $0<a<b$
\begin{equation}\label{eq::e12}
\hat{N}(u,a,b)=\sup_{a\le \rho\le b} \hat{N}(u,\rho)\quad \mbox{and}\quad \tilde{\omega}(f,a,b)=\sup_{a\le \rho\le b}\tilde{\omega}(f,\rho).
\end{equation}
Then we have the decay estimate below.
\begin{pro}\label{pro:decaybis}{\bf (Basic decay estimate)}\\
Given $p\in(1,+\infty)$, there exist constants $C_0>0$ $\lambda,\mu\in (0,1)$
(depending on $p$ and the dimension $n$) such that for every function $u$ and $f$ satisfying 
(\ref{eq:100}) with the notation given in (\ref{eq::e12}),  there holds
\begin{equation}\label{eq::e10}
\forall r\in (0,1],\quad 
\hat{N}(u,\lambda^2r,\lambda r)< \mu \ \hat{N}(u,\lambda r,r)
\quad \textup{or}\quad
\hat{N}(u,\lambda^2 r, \lambda r)< C_0 \ \tilde{\omega}(f,\lambda^2 r,r).
\end{equation}
\end{pro}
In order to prove this proposition, we will need the following result 
whose proof is postponed to subsection  \ref{sec:eqthm2}.

\begin{lem}\label{lem:Mestbis}{\bf (Estimates of $\hat{N}$ in larger balls)}\\
Let $u$ be a solution of 
$$\Delta u -u_t = f \quad \mbox{in}\quad B_R$$
for $R>2$, and
$$
\hat{N}(u,1)=||u-P_1||_{L^p(Q_1^-)}, 
$$
with $P_1\in {\P}_{2,1}$. Then for any $\rho \in [1,R/2]$, we have 
\begin{equation}\label{eq::e14}
\left(\frac{1}{\rho^{n+2+2p}}\int_{Q_\rho^-}|u-P_1|^p\d x\d t\right)^\frac{1}{p}\leq
  C_1\int_1^{2\rho}\frac{\hat{N}(u,s)+ \tilde{\omega}(f,s)}{s}\d s.
\end{equation}
\end{lem}

\noindent {\bf Proof of Proposition \ref{pro:decaybis}}\\ 
The proof is done by contradiction. If this is not true, we can
  find sequences $C_k\to \infty$, $r_k\in (0,1]$, $\lambda_k\to 0$ and $\mu_k\to 1$ such
  that (\ref{eq::e10}) fails with the corresponding functions
  $u_k$ and $f_k$ satisfying (\ref{eq:100}). This means that
\begin{equation}\label{eq::e16}
\left\{\begin{array}{l}
\hat N(u_k,\lambda_k^2r_k,\lambda_k r_k)\geq \mu_k \hat N (u_k,\lambda_k r_k,r_k),\\
\\
\hat  N(u_k,\lambda_k^2 r_k,\lambda_k r_k)\geq C_k\tilde \omega (f_k,\lambda_k^2 r_k, r_k).
\end{array}\right.
\end{equation}

\noindent {\bf Step 1: Construction of sequences and a priori estimates}\\
Let us consider a (not necessarily unique) $\rho_k \in [\lambda_k^2 r_k,\lambda_k r_k]$ so that
$$
\hat N(u_k,\lambda_k^2 r_k,\lambda_k r_k)=\hat N(u_k,\rho_k)=:\e_k.
$$
Moreover, define the rescaled functions
$$
v_k(x,t)=\frac{u_k(\rho_k x,\rho^2_k t)}{ \rho_k^2}
$$
and
\begin{equation}\label{eq:wk}
w_k(x,t)=\frac{u_k(\rho_k x,\rho^2_k t)-P_k(\rho_kx,\rho_k^2 t)}{\e_k \rho_k^2},
\end{equation}
where $P_k\in  {\mathcal P}_{2,c_{\rho_k}}$ is one polynomial realizing the infimum defining
$\hat N(u_k,\cdot)$ at the level $\rho_k$. Now, we wish to pass to the limit, but first we need to control the sequence $w_k$. By definition
\begin{equation}\label{eq:N1-0}
\inf_{P\in \P_{2}}\left(\int_{Q_1^-}\left|w_k-P\right|^p\d x \d t\right)^\frac{1}{p}=1.
\end{equation}
In addition, since
$$
\hat N(u_k,1)=\left(\frac{1}{\rho_k^{n+2+2p}}\int_{Q_{\rho_k}^-} |u_k-P_k|^p\d x\d t\right)^\frac{1}{p}, 
$$
we also have for $s\in (1,\frac{r_k}{2\rho_k})$ (applying Lemma \ref{lem:Mestbis} on $v_k$)
$$\begin{array}{ll}
\displaystyle \left(\frac{1}{s^{n+2+2p}}\int_{Q_s^-} |w_k|^p\d x\d
  t\right)^\frac{1}{p}&\displaystyle =\frac{1}{\e_k}\left(\frac{1}{(s\rho_k)^{n+2+2p}}\int_{Q_{s\rho_k}^-} |u_k-P_k|^p\d x\d t\right)^\frac{1}{p}\\
\\
&\displaystyle \le
\frac{C_1}{\e_k}\int_1^{2s}\frac{\hat N(u_k,\tau\rho_k) + \tilde{\omega}(f_k,\tau\rho_k)}{\tau}\d \tau\\
\\
&\displaystyle \leq \frac{C_1}{\e_k}\int_1^{2s}\frac{\hat N(u_k,\lambda_k^2 r_k,r_k)
+\tilde{\omega}(f_k,\lambda_k^2 r_k,r_k)}{\tau}\d \tau .
\end{array}$$
Notice that from (\ref{eq::e16}) we deduce
$$\hat N(u_k,\lambda_k^2 r_k,r_k) = \max(\hat N(u_k,\lambda_k^2 r_k,\lambda_kr_k), \hat N(u_k,\lambda_k r_k,r_k))\le \frac{\varepsilon_k}{\mu_k}$$
and
$$\tilde{\omega}(f_k,\lambda_k^2 r_k,r_k)\le \frac{\varepsilon_k}{C_k}.$$
This implies for $s\in (1,\frac{r_k}{2\rho_k})$ and some constant $C_2>0$
\begin{equation}\label{eq:subquad0}
\displaystyle \left(\frac{1}{s^{n+2+2p}}\int_{Q_s^-} |w_k|^p\d x\d
  t\right)^\frac{1}{p}\le C_2 \ln 2s.
\end{equation}
Furthermore, one can easily check that for  $H=\Delta -\partial_t$ and $s\in (1,\frac{r_k}{2\rho_k})$
we have
\begin{align*}
\left(\frac{1}{|Q_s^-|}\int_{Q_s^-}|H w_k|^p\right)^\frac1p
\leq \frac{1}{\e_k}\tilde \omega (f_k,s\rho_k)
\le\frac{1}{\e_k}\tilde \omega(f_k,\lambda_k^2r_k,r_k)
\leq \frac{1}{C_k}\to 0.
\end{align*}

\noindent {\bf Step 2: Identifying the limit and contradiction}\\
From \eqref{eq:subquad0} and the interior parabolic estimate (Theorem \ref{th:interior}), 
it follows that there is a subsequence again labeled $w_k$, converging in $L^p_\textup{loc}(\R^n\times \R^-)$ to a caloric 
function $w_0$. By passing to the limit in \eqref{eq:N1-0} we get
\begin{equation}\label{eq:N1limit}
\inf_{P\in \P_2}\left(\int_{Q_1^-}\left|w_0-P\right|^p\d x \d t\right)^\frac{1}{p}=1.
\end{equation}
Similarly, passing to the limit in \eqref{eq:subquad0} yields for all $s\ge 1$
$$
\left(\frac{1}{s^{n+2+2p}}\int_{Q_s^-} |w_0|^p\d x\d
  t\right)^\frac{1}{p}\leq C_2\ln 2s.
$$
Hence, $w_0$ is a caloric function in $\R^n\times \R^-$ that grows at most quadratically in space and linearly in time (up a logarithmic correction). This implies that $w_0$ is a caloric polynomial of degree at most two in space and one in time, i.e. $w_0\in\tilde  \P_2$. This clearly contradicts \eqref{eq:N1limit}.
This ends the proof of Proposition \ref{pro:decaybis}.\\

\subsection{Some routine results}\label{sec:eqthm2}

\noindent {\bf Proof of Lemma \ref{lem:Mestbis}}\\ 
The proof of this lemma is similar to the proof of Lemma 2.9 in \cite{Mon09}.\\
\noindent {\bf Step 1: Statement of (\ref{eq::e13})}\\
On the one hand, we use the fact that there exists a constant $C_2>0$ such that for any $P\in \tilde {\P}_2$ and for any $r\geq 1$
there holds
$$\left(\frac{1}{r^{n+2+2p}}\int_{Q_r^-} |P|^p\right)^{\frac1p}\le C_2 \left(\int_{Q_1^-} |P|^p\right)^{\frac1p}.$$
Following the proof of  Lemma 2.9 in \cite{Mon09}, 
we consider a dyadic decomposition of the cylinder $Q_\rho^-$,
and estimate  the quantities in each sub-cylinder.
More precisely, we get for $1\le \rho = 2^k r$ with $r\in [1/2,1)$ that
\begin{equation}\label{eq::e13}
\left(\frac{1}{\rho^{n+2+2p}}\int_{Q_\rho^-}|u-P_1|^p\d x\d t\right)^\frac{1}{p}\leq
  C_3 \left(\hat{N}(u,1)+\sum_{j=0}^k \hat{N}(u,2^j r)\right).
\end{equation}

\noindent {\bf Step 2: Proof of estimate (\ref{eq::e20})}\\ 
On the other hand, for any $\gamma>1$, we also notice that for $\alpha\in [1,\gamma]$, we have for any $r>0$ 
$$\hat{N}(u,\alpha r) \le \gamma^{2+\frac{n+2}{p}}\hat{N}(u,\gamma r) 
+ |c_{\alpha r}-c_{\gamma r}|\left(\int_{Q_1^-}|P_*|^p\right)^{\frac1p},$$
and
$$\left\{\begin{array}{l}
\displaystyle |c_{\alpha r}-c_{\gamma r}|
=\left(\frac{1}{|Q_{\alpha r}^-|}\int_{Q_{\alpha r}^-}|c_{\alpha r}-c_{\gamma r}|^p\right)^{\frac1p}
\le \tilde{\omega}(f,\alpha r) + \gamma^{\frac{n+2}{p}}\tilde{\omega}(f,\gamma r),\\
\\
\tilde{\omega}(f,\alpha r)\le \gamma^{\frac{n+2}{p}} \tilde{\omega}(f,\gamma r).
\end{array}\right.$$
Therefore, for any $\gamma>1$, there exists a constant $C_\gamma>0$ such that
\begin{equation}\label{eq::e20}
\forall \alpha\in [1,\gamma],\quad \left\{\begin{array}{l}
\hat{N}(u,\alpha r) \le C_\gamma \left(\hat{N}(u,\gamma r) +\tilde{\omega}(f,\gamma r)\right),\\
\\
\tilde{\omega}(f,\alpha r)\le   C_\gamma\tilde{\omega}(f,\gamma r).
\end{array}\right.
\end{equation}

\noindent {\bf Step 3: Conclusion}\\ 
Using (\ref{eq::e20}) with $\gamma=2$, we 
get the result (\ref{eq::e14}) with the integral on $\hat{N}$ 
replacing the sum in the right hand side of  (\ref{eq::e13}). This ends the proof of Lemma \ref{lem:Mestbis}.\\


\noindent Given $(u,f)$ and $\lambda\in (0,1)$, let us introduce the notation
\begin{equation}\label{eq::e18}
\underline{N}(r)=\hat N(u,\lambda r,r)\textup{ and } \underline{\omega}(r)=\tilde{\omega}(f,\lambda^2 r,r).
\end{equation}
Contrarily to what is done in \cite{Mon09}, 
the functions $\underline{N}$ and $\underline{\omega}$ are not necessarily monotone in $r$. 
Nevertheless, we have the following routine result (the analogue to Lemma 3.4 in \cite{Mon09}).

\begin{pro}{\bf (Dini estimate)}\label{pro::e17}\\
Let $\underline{N}:(0,1]\to [0,+\infty)$, $\underline{\omega}:(0,1]\to [0,+\infty)$
be two functions satisfying
\begin{equation}\label{eq::e19}
\forall r\in (0,1],\quad \underline{N}(\lambda r)\le \mu \underline{N}(r)\quad \mbox{or}\quad \underline{N}(\lambda r)\le \underline{C}_0\ \underline{\omega}(r)
\end{equation}
and
\begin{equation}\label{eq::e21}
\forall r\in (0,1],\quad \forall \alpha\in [\lambda,1],\quad \left\{\begin{array}{l}
\underline{N}(\alpha r) \le \underline{C}_0\left(\underline{N}(r) +\underline{\omega}(r)\right),\\
\\
\underline{\omega}(\alpha r)\le \underline{C}_0 \underline{\omega}(r),
\end{array}\right.
\end{equation}
for some constants $\underline{C}_0>0$, $\lambda,\mu\in (0,1)$ and assume that $\underline{\omega}$ is Dini.
Then there exists a constant $\underline{C}_0'>0$ depending only on $\underline{C}_0,\lambda,\mu >0$, such that for every $\rho\in (0,\lambda/2]$ and with $\alpha=\ln \mu/\ln \lambda$ there holds
$$\int_0^\rho \frac{\underline{N}(r)}{r}\ dr
\le \underline{C}_0' \left\{\underline{N}(1) \rho^\alpha + \int_0^\rho \frac{\underline{\omega}(r)}{r}\ dr + \rho^\alpha  \int_\rho^{1}\frac{\underline{\omega}(r)}{r^{1+\alpha}}\ dr\right\}.$$
\end{pro}

\begin{rem}
Notice that the quantities $\underline{N}$ and $\underline{\omega}$ defined in (\ref{eq::e18})
satisfy (\ref{eq::e19}) because of the basic estimate (Proposition \ref{pro:decaybis}) and 
 do also satisfy (\ref{eq::e21}) because of (\ref{eq::e20}) with $\gamma =1/\lambda$ (with $\underline{C}_0=\max(C_0,C_{1/\lambda})$).
\end{rem}

\noindent {\bf Proof of Proposition \ref{pro::e17}}\\
\noindent {\bf Step 1: Estimate on $\underline{N}(r)$}\\
We claim that we have for all $r\in (0,\lambda]$
\begin{equation}\label{eq::e23}
\underline{N}(r)\le \max \left(C_2 r^\alpha, \underline{C}_0 \ r^\alpha 
\sup_{\rho\in [r,\lambda]} \frac{\underline{\omega}(\rho)}{\rho^\alpha}\right).
\end{equation}
The proof is the same as Lemma 3.3 in \cite{Mon09} with $r_0=1$, except that we estimate for $r_1\in (\lambda,1]$
$$\underline{N}(r_1) \le \underline{C}_0(\underline{N}(1)+\underline{\omega}(1)).$$
This gives the new value to the constant 
\begin{equation}\label{eq::e24}
C_2= \lambda^{-\alpha}\underline{C}_0(\underline{N}(1)+\underline{\omega}(1)).
\end{equation}
Here we have replaced the lack of monotonicity of $\underline{N}$ 
by the first line of assumption (\ref{eq::e21}).\\
\noindent {\bf Step 2: Estimate on $\displaystyle \sup_{\rho\in [r,\lambda]} \frac{\underline{\omega}(\rho)}{\rho^\alpha}$
and conclusion}\\
We follow the proof of Lemma 3.4 in \cite{Mon09}. For some $\rho_0\in [r,\lambda]$ we have
$$\begin{array}{lll}
\displaystyle \sup_{\rho\in [r,\lambda]} \frac{\underline{\omega}(\rho)}{\rho^\alpha}
&\displaystyle  =   \frac{\underline{\omega}(\rho_0)}{\rho_0^\alpha} \\
\\
& \displaystyle \le \frac{1}{\rho_0^\alpha} \frac{1}{t\rho_0}\int_{\rho_0}^{\rho_0+t\rho_0} C_0 
\underline{\omega}(\rho)\ d\rho
&  \displaystyle \quad \mbox{with}\quad t=\frac{1-\lambda}{\lambda} >0\\
\\
& \displaystyle \le \frac{\underline{C}_0}{t\lambda^{1+\alpha}}\int_{\rho_0}^{\frac{\rho_0}{\lambda}} \frac{\underline{\omega}(\rho)}{\rho^{1+\alpha}}\ d\rho& \\
\\
& \displaystyle \le C_3 \int_r^1 \frac{\underline{\omega}(\rho)}{\rho^{1+\alpha}}\ d\rho&  
\displaystyle \quad \mbox{with}\quad C_3=\frac{\underline{C}_0}{(1-\lambda)\lambda^\alpha}>0 ,
\end{array}$$
where in the second line we have used the second line of assumption (\ref{eq::e21})
(because of the lack of monotonicity of $\underline{\omega}$).
The remaining part of the proof of Lemma 3.4 in \cite{Mon09} is unchanged and then implies the result.
This ends the proof of Proposition \ref{pro::e17}.\\

\subsection{Proof of Theorem \ref{th:1-bis}}\label{sec:eqthm3}

\noindent {\bf Proof of Theorem \ref{th:1-bis}}\\ 
\noindent {\bf Proof of i)}\\
Using definition (\ref{eq::e18}) of $\underline{N}$ and $\underline{\omega}$, and estimate (\ref{eq::e23})
with the constant $C_2$ given in (\ref{eq::e24}), we deduce that for $r\in (0,\lambda]$
$$\underline{N}(r)\le \underline{C}_0 
\left(\underline{N}(1)+ \sup_{\rho\in (0,1]} \underline{\omega}(\rho)\right).$$
From (\ref{eq::e20}) with $\gamma=1/\lambda$, we deduce that for all $r\in (0,1]$
$$\hat{N}(u,r) \le C\left\{\hat{N}(u,1) + \sup_{\rho\in (0,1]} \tilde{\omega}(f,\rho)\right\}$$
which implies (\ref{eq::e25}) because we always have
$$\tilde{N}(u,r)\le \hat{N}(u,r)$$
and 
$$
 \hat N(u,1)\leq C\left(\|u\|_{L^p(Q_1^-)}+\|f\|_{L^p(Q_1^-)}\right).
$$
\noindent{\bf Proof of ii)}\\
This follows from (\ref{eq::e23}).\\
\noindent {\bf Proof of iii)}\\
\noindent {\bf Step 1: Dini property}\\
Proposition \ref{pro::e17} implies that $\underline{N}$ is Dini if $\underline{\omega}$ is Dini.
Then we deduce that $\hat{N}(u,r)$ (and then $\tilde{N}(u,r)$) is Dini, if $\tilde{\omega}(f,r)$ is Dini.\\
\noindent {\bf Step 2: Estimate (\ref{eq::e26})}\\
The proof of Lemma 3.5 in \cite{Mon09} is straightforward to adapt to our case.
Using our Lemma \ref{lem:Mestbis} instead of Lemma 2.9 in \cite{Mon09}, 
this shows that there exists a polynomial $P_0\in \tilde {\P}_2$ such that for all $\rho\in (0,r_*]$ with $r_*=\lambda/2$
$$\left(\frac{1}{\rho^{n+2+2p}}\int_{Q_\rho^-}|u-P_0|^p\d x\d t\right)^\frac{1}{p}\leq
C_1\int_1^{2\rho}\frac{\hat{N}(u,s)+\tilde{\omega}(f,s)}{s}\d s.$$
We deduce (\ref{eq::e26}) using Proposition \ref{pro::e17} joint with (\ref{eq::e20}).\\
\noindent {\bf Step 3: $P_0$ is caloric}\\
From (\ref{eq::e26}) and the interior estimates (Theorem \ref{th:interior}) we also deduce that
$$u^\varepsilon(x,t)=\frac{(u-P_0)(\varepsilon x,\varepsilon^2 t)}{\varepsilon^2}$$
converges, as $\e\to 0$, to a function $v\equiv 0$, which is a solution of
$$Hv = f(0)- H(P_0) \quad \mbox{with}\quad H =\Delta -\partial_t.$$
Since $f(0)=0$, this shows that $P_0$ is caloric.\\
\noindent {\bf Step 4: Bound on the coefficients of $P_0$}\\
We simply apply (\ref{eq::e26}) for $r=\lambda/2$ and this implies the bound on the coefficients of $P_0$.\\
This ends the proof of the theorem.

\section{Growth estimates for the obstacle problem}\label{sec:growth}

In this section we will prove some growth estimates of solutions to
\eqref{eq:obstacle}. Some of the results are of independent interest
while some are needed in the sequel.

\begin{pro}\label{pro:quadgrowth} {\bf (Quadratic growth in mean)}\\ 
Let $u$ be a solution of
  \eqref{eq:obstacle} with $p\in ((n+2)/2,\infty)$. Then there are positive constants $r_0$ and 
$$
 C_1=C\left(1+\sigma(1)\right)\|u\|_{L^p(Q_1^-)},
$$ such that 
$$
\left(\frac{1}{|Q_r^-|}\int_{Q_r^-}|u|^p\right)^\frac1p\leq C_1r^2, 
$$
whenever $r<r_0$.
\end{pro}
\noindent {\bf Proof of Proposition \ref{pro:quadgrowth}}\\
  Define for $r\in (0,1]$
$$
S_r(u)=\left(\frac{1}{|Q_r^-|}\int_{Q_r^-}|u|^p\right)^{\frac1p}.
$$
By iteration it is sufficient to prove that there exists $C>0$ and $r_0\in (0,1]$ such that for all solutions $u$ of \eqref{eq:obstacle}, for all $r\leq r_0$ either
$$ S_r(u)\le C r^2\left(1+ \sigma(f,1)\right)$$
or
$$
S_r(u)\le 4^{-k} S_{2^k r}(u) \quad \mbox{for some $k\in \N$ such that}\quad 2^k r \le 1.$$
In order to
prove that this holds we argue by contradiction. If this does not hold,
there are sequences $u_j$, $f_j$ and $r_j\to 0$,$C_j\to \infty$ such that 
$$ S_{r_j}(u_j)\geq C_j r_j^2\left(1+ \sigma(f_j,1)\right)$$
and
$$ S_{r_j}(u_j)\geq 4^{-k}S_{2^kr_j}(u_j) 
\quad \mbox{for all}\quad  k\in \N\quad \mbox{such that} \quad 2^kr_j\le 1.
$$
Define the rescaled
functions
$$
v_j(x,t)=\frac{u_j(r_jx,r_j^2t)}{S_{r_j}(u_j)}.
$$
Then
\begin{enumerate}
\item $\displaystyle Hv_j=\frac{r_j^2}{S_{r_j}(u_j)}\ f_j(r_jx,r_j^2t)\chi_{\left\{v_j>0\right\}}$ in $Q_{\frac{1}{r_j}}^-$,
\item $S_1(v_j)=1$,
\item $S_{2^k}(v_j)\leq 4^k$ for all $k\in \N$ with $2^kr_j\le 1$,
\item $v_j(0)=0$,
\item $v_j\geq 0$.
\end{enumerate}
Observing that for $\rho\le 1/r_j$ we have
$$
\left(\frac{1}{|Q_\rho^-|}\int_{Q_\rho^-}|f_j(r_jx,r_j^2t)|^p\right)^{\frac1p}
=\left(\frac{1}{|Q_{\rho r_j}^-|}\int_{Q_{\rho r_j}^-}|f_j(x,t)|^p\right)^{\frac1p}
\le 1 + \sigma(f_j,1),
$$
we see that $Hv_j$ is bounded in $L^p(Q_\rho^-)$ for every
$\rho$ and converges to 0. Therefore, by interior parabolic estimates and the Sobolev
embedding (Theorem \ref{th:interior} and Theorem \ref{th:sobolev})
there exists a subsequence, again labelled $v_j$, such that $v_j\to v_0$
locally in $C^\alpha$ and locally weakly in
$W^{2,1}_p$, where $v_0$ satisfies
\begin{enumerate}
\item $Hv_0=0$ in $\R^n\times\R^-$,
\item $S_1(v_0)=1$,
\item $v_0(0)=0$,
\item $v_0\geq 0$.
\end{enumerate}
This contradicts the strong maximum principle for caloric functions (see Theorem 11 on page 375 in \cite{evans2009partial}) and ends the proof of the proposition.\\

Using Proposition \ref{pro:quadgrowth} we prove the following corollary that implies Proposition \ref{pro:sing1}.

\begin{cor}\label{cor:quadgrowth}{\bf (Quadratic growth in the $\sup$-norm)}\\
Under the same assumptions as in
  Proposition \ref{pro:quadgrowth}, there is a constant $C_2>0$ such that there holds
$$
\sup_{Q_r^-}|u|\leq C_2C_1r^2 \quad \mbox{for all}\quad r < r_0/2,
$$
where $r_0$ and $C_1$ are defined in Proposition \ref{pro:quadgrowth}.
\end{cor}
\noindent {\bf Proof of Corollary \ref{cor:quadgrowth}}\\ 
Define
$$
v_r(x,t)=\frac{u(rx,r^2t)}{r^2}.
$$
Then Proposition \ref{pro:quadgrowth} implies that for $r<r_0$,
$$
\left(\frac{1}{|Q_1^-|}\int_{Q_1^-}|v_r|^p\right)^{\frac1p}\leq C_1.
$$
Moreover, there holds
$$
\displaystyle \left(\frac{1}{|Q_1^-|}\int_{Q_1^-}|Hv_r|^p\right)^{\frac1p}\leq 1 + \sigma(1).
$$  
Therefore, by interior estimates (Theorem \ref{th:interior})
$$
||v_r||_{W^{2,1}_p}(Q_\frac12^-)\leq C_2' C_1,
$$
and thus, by the Sobolev embedding (Theorem \ref{th:sobolev})
$$
\sup_{Q_\frac12^-}|v_r|\leq C_2C_1.
$$
Scaling back to $u$ yields the desired result.


\section{Non-degeneracy}\label{sec:nondeg}

In this section we prove that solutions of \eqref{eq:obstacle} cannot
decay too fast close to the origin, and the rate at which this can
happen, naturally depends on $f$.


\begin{pro}{\bf (Non-degeneracy)}\label{pro:nondeg} \\
Let $p\in
  ((n+2)/2,\infty)$. In addition, assume that $u$ solves
$$\left\{
\begin{array}{l}
\left.\begin{array}{l}
Hu=f\chi_{\{u>0\}}\\
u\geq 0
\end{array}\right| \textup{ in $Q_R^-$,}\\
\ f(0)=1.
\end{array}\right.$$
Then there exists a constant $C_0>0$ such that
$$
\sup_{\dd_p Q_d^-(x_0,t_0)}u\geq \frac{d^2}{2n+1},
$$
whenever $Q_d^-(x_0,t_0)\subset Q_R^-$ and $u(x_0,t_0)>2\lambda$, 
where
$$
\lambda=C_0R^\frac{n+2}{p}d^{2-\frac{n+2}{p}}\left(\frac{1}{|Q_R^-|}\int_{Q_R^-} |f(x,t)-f(0)|^p\right)^\frac{1}{p}.
$$
\end{pro}

\noindent {\bf Proof of Proposition \ref{pro:nondeg}}\\ 
Let $v$ be the solution of 
$$\left\{\begin{array}{lll}
Hv=f-f(0)  &\quad\mbox{in}\quad & Q_d^-(x_0,t_0),\\
v= 0 & \quad \mbox{on}\quad  &\dd_p Q_d^-(x_0,t_0).
\end{array}\right.$$
Defining
$$
v_d(x,t)=\frac{v(dx+x_0,d^2t+t_0)}{d^2}
$$
and $f_d(x,t)=f(dx+x_0,d^2t+t_0)$,
we see that
$$\left\{\begin{array}{lll}
Hv_d=f_d-f(0)  &\quad\mbox{in}\quad & Q_1^-,\\
v_d= 0 & \quad \mbox{on}\quad  &\dd_p Q_1^-.
\end{array}\right.$$
The classical parabolic estimates (Theorem \ref{th:para}) imply
$$
||v_d||_{W^{2,1}_p(Q_1^-)}\leq C||f_d-f(0)||_{L^p(Q_1^-)}.
$$
Applying the Sobolev embedding (Theorem \ref{th:sobolev}) and scaling back to $v$, we deduce
$$
\sup_{Q_d^-(x_0,t_0)}|v|\leq d^2C_0 \left(\frac{1}{d^{n+2}|Q_1^-|}\int_{Q_d^-(x_0,t_0)} |f(x,t)-f(0)|^p\right)^\frac{1}{p} \le \lambda,
$$
with $C_0=C_*C|Q_1^-|^{\frac1p}$. Let 
$$
w=u-v-f(0)\left(\frac{(x-x_0)^2-(t-t_0)}{2n+1}\right).
$$
Then
$$Hw= 0\quad \mbox{in}\quad  Q_d^-(x_0,t_0)\cap\{u>0\}.$$
Moreover, $w(x_0,t_0)=u(x_0,t_0)-v(x_0,t_0)>\lambda$. Therefore, by the
maximum principle, $w$ attains its
positive maximum $>\lambda$ on $\dd_p(Q_d(x_0,t_0)\cap
\{u>0\})$. Whenever $u=0$ we have $w\le \lambda$. Thus, the maximum is
attained on $\{u>0\}\cap \dd_pQ_d(x_0,t_0)$ and we have
$$
\lambda<\sup_{\{u>0\}\cap \dd_pQ_d(x_0,t_0)}w\leq \sup_{\{u>0\}\cap \dd_pQ_d(x_0,t_0)}u-\frac{d^2}{2n+1} +\lambda.
$$
The result follows.

A corollary from this non-degeneracy follows below.
\begin{cor}\label{cor:weaknondeg}{\bf (Weak non-degeneracy)}\\
Let $p\in ((n+2)/2,\infty)$. Assume that $u_m$ and $f_m$ verify
$$\left\{
\begin{array}{l}
\left.\begin{array}{l}
Hu_m=f_m\chi_{\{u_m>0\}}\\
u_m\geq 0
\end{array}\right| \quad \mbox{in}\quad Q_R^-,\\
\\
\ f_m(0)=1,\\
\left.\begin{array}{l}
\displaystyle \tau_m=\left(\frac{1}{|Q_R^-|}\int_{Q_R^-}|f_m(x,t)-f_m(0)|^p\right)^\frac{1}{p}\to
    0\\
u_m\to u_\infty \quad \mbox{in}\quad  L^\infty_{\textup{loc}}(Q_R^-)
\end{array}\right| \quad \mbox{as}\quad  m\to \infty.
\end{array}\right.
$$
Then for any compact $K\subset \{u_\infty=0\}^\circ\cap Q_R^-$ (where $\{u_\infty=0\}^\circ$ denotes the interior of the set $\{u_\infty=0\}$), there exists a
constant $C>0$ (independent of $\tau_m$) such that 
$$
u_m\leq C\tau_m\textup{ in $K$}.
$$
\end{cor}
\noindent {\bf Proof of Corollary \ref{cor:weaknondeg}}\\ 
We argue  by
  contradiction. Choose $d$ such that 
$$
\left(\left(\bigcup_{P\in K}Q_d^-(P)\right)\bigcap \{t\leq 0\}\right)\quad \subset\subset \quad  Q_R^-\cap \{u_\infty=0\}.
$$
Suppose that $u_m(P_m)\geq C_m\tau_m$ for $P_m\in K$ and $C_m\to
\infty$. Clearly, for $m$ large enough we have
$$
\lambda_m=C_0R^\frac{n+2}{p}d^{2-\frac{n+2}{p}}\left(\frac{1}{|Q_R^-|}\int_{Q_R^-} |f_m(x,t)-f_m(0)|^p\right)^\frac{1}{p} < \frac{C_m\tau_m}{2},
$$
and thus
$$
u_m(P_m)>2\lambda_m.
$$
Then, from Proposition \ref{pro:nondeg}, we know that
$$
\sup_{\overline{Q_d^-(P_m)}}u_m\geq \frac{d^2}{2n+1}, 
$$
which implies
$$
\sup_{\overline{Q_d^-(P_\infty)}}u_\infty\geq \frac{d^2}{2n+1},
$$
where $P_m\to P_\infty\in K$. This is a contradiction.

\section{A compactness result}\label{sec:comp}

The main result of this section will be Corollary \ref{cor:conv} which shows the compactness in $L^p$ of certain sequences.
This result will be applied in the next section.

\begin{lem}{\bf (Cacciopoli type estimate)}\label{lem:cac1} \\
Let $u$ be a
  solution of
\begin{equation}\label{eq:eq1}
\left\{\begin{array}{l}
\left.\begin{array}{l}
Hu=f\chi_{\{u>0\}}\\
u\geq 0
\end{array}\right|\quad \mbox{in}\quad Q_1^-,\\
\ f(0)\ge 0,
\end{array}\right.
\end{equation}
and $P$ be a solution of (\ref{eq:eq1}) with $f$ replaced by the constant function $f(0)$.
Furthermore, set $w=u-P$ and
$W=w|w|^{\frac{p}{2}-1}$ for $p\in (1,+\infty)$. Then for any
$\eta\in C_0^\infty(\R^{n+1})$ such that $\supp\eta\subset Q_r$ with
$0<r\leq 1$, we have
\begin{equation}\label{eq:r100}
\begin{array}{l}
\displaystyle \int_{Q_1^-}\frac{2(p-1)}{p^2}|\nabla
W|^2\eta^2+f(0)(\chi_{\{u>0\}}-\chi_{\{P>0\}})|w|^{p-2}w\eta^2+\int_{B_1\times\{0\}}\eta^2\frac{W^2}{p}\\
\\
\displaystyle \leq
\int_{Q_1^-}\frac{2}{p-1}W^2|\nabla\eta|^2+\frac{2}{p}|\eta\eta_t|W^2+\omega(r)|Q_r^-|\left(\frac{1}{|Q_r^-|}\int_{Q_r^-}\eta^{2p'}W^2\right)^\frac{1}{p'}.
\end{array}
\end{equation}
\end{lem}

\noindent {\bf Proof of Lemma \ref{lem:cac1}}\\ 
Recall that $w=u-P$ solves the equation
$$Hw=f\chi_{\{u>0\}}-f(0)\chi_{\{P>0\}}\quad \mbox{in}\quad  Q_1^-.$$
Multiplying this equation with $\eta^2w|w|^{p-2}$ we get
$$\int_{Q_1^-}(Hw)\eta^2w|w|^{p-2}=\int_{Q_1^-} \left(f(0)(\chi_{\{u>0\}}-\chi_{\{P>0\}})w|w|^{p-2}\eta^2 +(f-f(0))\chi_{\{u>0\}}w|w|^{p-2}\eta^2\right).$$
Furthermore, 
$$\int_{Q_1^-}(Hw)\eta^2w|w|^{p-2}
=\int_{Q_1^-} \eta^2w|w|^{p-2}\lap
w+2\eta\eta_t\frac{|w|^p}{p}-\int_{Q_1^-\cap
  \{t=0\}}\eta^2\frac{|w|^p}{p},$$
and
$$\int_{Q_1^-}-\eta^2w|w|^{p-2}\lap w
=\int_{Q_1^-}\frac{4(p-1)}{p^2}\eta^2|\nabla
W|^2+\frac{4}{p}\eta W\nabla \eta\nabla W.$$
Therefore, with $\lambda=\frac{4(p-1)}{p^2}$ we have
$$\begin{array}{ll}
&\displaystyle \int_{Q_1^-}\lambda\eta^2|\nabla W|^2+f(0)(\chi_{\{u>0\}}-\chi_{\{P>0\}})|w|^{p-2}w\eta^2+\int_{Q_1^-\cap
  \{t=0\}}\eta^2\frac{W^2}{p}\\
\\
\leq&\displaystyle \int_{Q_1^-} -\frac{4}{p}\eta W\nabla
\eta\nabla W+2\eta\eta_t\frac{W^2}{p}+\eta^2|f-f(0)|W^\frac{2(p-1)}{p}\\
\\
\leq& \displaystyle \int_{Q_1^-}\frac{1}{2}\left(\lambda\eta^2|\nabla W|^2+\frac{16}{p^2}\lambda^{-1}W^2|\nabla\eta|^2\right)+2\eta\eta_t\frac{W^2}{p}+\omega(r)|Q_r^-|\left(\frac{1}{|Q_r^-|}\int_{Q_r^-}\eta^{2p'}W^2\right)^\frac{1}{p'}.
\end{array}$$
This implies
\begin{align*}
  \int_{Q_1^-}\frac{1}{2}\lambda\eta^2|\nabla W|^2+f(0)(\chi_{\{u>0\}}-\chi_{\{P>0\}})|w|^{p-2}w\eta^2+\int_{Q_1^-\cap
  \{t=0\}}\eta^2\frac{W^2}{p}\\
\leq \int_{Q_1^-}\frac{8}{p^2}\lambda^{-1}W^2|\nabla\eta|^2+2|\eta\eta_t|\frac{W^2}{p}+\omega(r)|Q_r^-|\left(\frac{1}{|Q_r^-|}\int_{Q_r^-}\eta^{2p'}W^2\right)^\frac{1}{p'},
\end{align*}
which is the desired inequality.

\begin{lem}{\bf ($L_t^\infty L_x^2$-estimates for $W$)}\label{lem:Linfty} \\
Under the same
  assumptions as in Lemma \ref{lem:cac1}, we have for any $t_0\in
  [-\frac{1}{4},0]$
$$\int_{Q_\frac{1}{2}^-\cap \{t=t_0\}}W^2\leq C\left(\int_{Q_1^-}W^2+(\omega(1))^p\right).$$
\end{lem}

\noindent {\bf Proof of Lemma \ref{lem:Linfty}}\\ 
Let $w^{t_0}(x,t)=w(x,t-t_0)$ with $u^{t_0}(x,t)=u(x,t-t_0)$ and $P^{t_0}(x,t)=P(x,t-t_0)$. Then we have
$$Hw^{t_0}=f^{t_0}\chi_{\{u^{t_0}>0\}}-f(0)\chi_{\{P^{t_0}>0\}}\quad \mbox{in}\quad  Q_{\sqrt{1-t_0^2}}^-.$$
Put $W^{t_0}=W(x,t-t_0)$ and take $\eta\in
  C_0^\infty(\R^{n+1})$ such that  $\supp\eta\subset Q_r$ with $r>1/2$
and $1/4+r^2<1$ (so that $Q_r^-\subset Q_{\sqrt{1-t_0^2}}^-$), and $\eta=1$ on $B_\frac{1}{2}\times\{t=0\}$. We now apply the proof of Lemma \ref{lem:cac1} to $W^{t_0}$ together with Young's
inequality applied to the last term of (\ref{eq:r100}). This gives
$$
\int_{Q_1^-\cap \{t=0\}}(W^{t_0})^2\eta^2\leq C\left(\int_{Q_r^-}(W^{t_0})^2+\int_{Q_r^-}|f^{t_0}-f(0)|^p\right).
$$
Since
$$
\int_{Q_r^-}(W^{t_0})^2\leq \int_{Q_1^-}W^2, 
$$
$$
\int_{Q_r^-}|f^{t_0}-f(0)|^p\leq \int_{Q_1^-}|f-f(0)|^p=|Q_1^-|(\omega(1))^p
$$
and
$$
\int_{Q_\frac{1}{2}^-\cap \{t=t_0\}}W^2\leq \int_{Q_\frac{1}{2}^-\cap
  \{t=0\}}(W^{t_0})^2\eta^2, 
$$
this implies the result.

\begin{cor}\label{cor:w1q}{\bf ($L^q_tW_x^{1,q}$-estimates for $w$)}\\
Under the same assumptions as in Lemma \ref{lem:cac1}, there is $q\in (1,p]$ such that for any $0<r<1$, we have
$$\int_{Q_r^-}|\nabla w|^q\leq C(p,r,\omega(1),||w||_{L^p(Q_1^-)}).$$
\end{cor}

\noindent {\bf Proof of Corollary \ref{cor:w1q}}\\ 
We divide the proof into two cases, depending on whether
  $p>2$ or not. In order to clarify the dependence of $\omega$ on $p$, we write it $\omega_p$.

\noindent {\bf Case 1:} $p\geq 2$. In this case we can simply apply Lemma \ref{lem:cac1}
with $p=2$ to obtain for $0<r<1$:
$$\int_{Q_r^-}|\nabla w|^2\leq C(r,\omega_2(1),||w||_{L^2(Q_1^-)})\leq C(r,\omega_p(1),||w||_{L^p(Q_1^-)})$$

\noindent {\bf Case 2:} $p<2$. We compute for $0<r<1$:
$$\begin{array}{ll}
\displaystyle \int_{Q_r^-}|\nabla w|^q & \displaystyle = \int_{Q_r^-}|\nabla w|^q|w|^{q(p/2-1)}|w|^{-q(p/2-1)}\\
\\
&\displaystyle \leq\left(\int_{Q_r^-}|\nabla w|^2 w^{p-2}\right)^{q/2}\left(\int_{Q_r^-}(|w|^{-q(p/2-1)})^{(2/q)'}\right)^{1/(2/q)'}\\
\\
&\displaystyle \leq C(r,\omega_p(1),||w||_{L^p(Q_r^-)}) 
\end{array}$$
The first factor can be estimated using Lemma \ref{lem:cac1} and the
second one by $C||w||_{L^p(Q_r^-)}$. This is due to the fact that
$(2/q)'=(1-q/2)^{-1}$ which implies that the exponent of $w$ is
$-q(p/2-1)/(1-q/2)=-q(p-2)/(2-q)$. We realize that for $q=1$ this equals
$2-p<p$, and hence if we take $q>1$ small enough the exponent will be
less than $p$.

\begin{lem}\label{lem:L1}{\bf (Partial $L^1$-estimates of the right hand side)}\\ 
Under the same assumptions as in Lemma \ref{lem:cac1}, for any $0<r<1$ there holds
\begin{equation}\label{eq:r101}
\int_{Q_r^-}f(0)|(\chi_{\{u>0\}}-\chi_{\{P>0\}})|\leq C(p,r,\omega(1),||w||_{L^p(Q_1^-)}).
\end{equation}
\end{lem}

\noindent {\bf Proof of Lemma \ref{lem:L1}}\\
For $\e\geq 0$ define
$$
\sgn_\e(w)=\left\{\begin{array}{lr}
1 & \textup{ if $w>\e$,}\\
\frac{w}{\e} & \textup{ if $|w|\leq\e$,}\\
-1 & \textup{ if $w<-\e$,}\\
\end{array}\right.
$$
and $h_\e=f(0)(\chi_{\{u>0\}}-\chi_{\{P>0\}})\sgn_\e(w)$. As before we
have
$$
Hw=(f-f(0))\chi_{\{u>0\}}+f(0)(\chi_{\{u>0\}}-\chi_{\{P>0\}}).
$$
Take $\eta\in C_0^\infty(Q_\rho^-)$ with $\eta=1$ on $Q_r^-$ for
$r<\rho<1$. Multiplying the equation by $\sgn_\e(w)\eta^2$ yields
$$
\int_{Q_1^-}h_\e\eta^2=-\int_{Q_1^-}(f-f(0))\chi_{\{u>0\}}\eta^2\sgn_\e(w)+\int_{Q_1^-}(Hw)\sgn_\e(w)\eta^2.
$$
We observe that $\sgn_\e(w)$ is the derivative of a convex function
$\beta_\e(w)\ge 0$. Therefore
$$
\int_{Q_1^-}(Hw)\sgn_\e(w)\eta^2=\int_{Q_1^-}\eta^2\sgn_\e(w)\lap
w+2\eta\eta_t\beta_\e(w)-\int_{Q_1^-\cap \{t=0\}}\eta^2\beta_\e(w), 
$$
and
$$
-\int_{Q_1^-}\eta^2\sgn_\e(w)\lap w=\int_{Q_1^-}|\nabla
w|^2\frac{1}{\e}\chi_{\{|w|\leq \e\}}\eta^2+2\eta\nabla \eta\nabla w\sgn_\e(w).
$$
Adding up, this gives
\begin{align*}
\int_{Q_1^-}h_\e\eta^2+\int_{Q_1^-\cap
  \{t=0\}}\eta^2\beta_\e(w)+\int_{Q_1^-}|\nabla
w|^2\frac{1}{\e}\chi_{\{|w|\leq \e\}}\eta^2\\\leq
\int_{Q_1^-}|f-f(0)|\eta^2+2|\eta\eta_t|\beta_\e(w)+2\eta|\nabla
\eta\nabla w|. 
\end{align*}
Observing that $|\beta_\e(w)-|w||\leq
\e/2$, we see
$$\int_{Q_r^-}h_\e\leq C\omega(1)+C(r)\left(||w||_{L^1(Q_1^-)}+\e + ||\nabla w||_{L^1(Q_\rho^-)}\right).$$
Hence, the Dominated convergence theorem implies $h_\e\to
h_0=f(0)\left(\chi_{\left\{u>0\right\}}-\chi_{\left\{P>0\right\}}\right)\sgn(w)$ in $L^1(Q_r^-)$. 
Thus, by Corollary  \ref{cor:w1q},
$$\int_{Q_r^-}h_0\leq C(p,r,\omega(1),||w||_{L^p(Q_1)}),$$
which gives (\ref{eq:r101}). This ends the proof of the lemma.\\

Below we state a result of Simon we will be using (see Theorem 6 on page 86 in \cite{Sim87}) and a small lemma that we will need.

\begin{theo}\label{th:aubin}{\bf (Compactness in Banach spaces)}\\
Let $X_0\subset X\subset X_1$ be Banach spaces such that $X_0$ is
compactly embedded in $X$ and $X$ is continuously embedded in
$X_1$. Moreover, assume that $u_k$ is a sequence of functions such that
for some $q>1$ 
$$
||u_k||_{L^q(I;X)}+||u_k||_{L^1(I;X_0)}+||\dd_t u_k||_{L^1(I;X_1)}\leq C, 
$$
where $I\subset \R $ is a compact interval. Then there is a subsequence
$u_{k_j}$ that converges in $L^p(I;X)$ for all $1\leq p<q$.
\end{theo}

\begin{lem}\label{lem:inclusion}{\bf (Inclusion in dual spaces)}\\
For any $r$ there holds
$$
\nabla L^s(B_r)\subset W^{-1,s}(B_r)
$$
and
$$
L^1(B_r)\subset W^{-1,s}(B_r), 
$$
whenever
$$
1<s<\frac{n}{n-1}.
$$
\end{lem}

\noindent {\bf Proof of Lemma \ref{lem:inclusion}}\\
We first we prove the inclusion
$$
\nabla L^s(B_r)\subset W^{-1,s}(B_r),
$$
where $W^{-1,s}(B_r)$ is by definition the dual of the space $W^{1,s'}_0(B_r)$ with $1/s + 1/s' =1$.
Let $f\in L^s(B_r)$ and $\phi\in C_0^\infty(B_r)$. Then, by H\"older's inequality
$$
|\<\nabla f,\phi\>|=|\<f,\nabla \phi\>|\leq ||f||_{L^s(B_r)}||\phi||_{W_0^{1,s'}(B_r)}.
$$
Since $C_0^\infty(B_r)$ is dense in $W^{1,s'}(B_r)$, this implies
$$
||\nabla f||_{W^{-1,s}(B_r)}=||\nabla f||_{(W_0^{1,s'}(B_r))^*}\leq ||f||_{L^s(B_r)}.
$$
Now we prove that $L^1(B_r)\subset W^{-1,s}(B_r)$ when $s<\frac{n}{n-1}$. In order
to do so, take $g\in L^1(B_r)$ and $\phi\in W_0^{1,s'}$. H\"older's inequality implies
$$
\left|\int_{B_r}g\phi\right|\leq ||g||_{L^1(B_r)}||\phi||_{L^\infty(B_r)}\leq
  ||g||_{L^1(B_r)}||\phi||_{W^{1,s'}(B_r)}, 
$$
whenever $s'>n$, which is equivalent to $s<\frac{n}{n-1}$.
This ends the proof of the lemma.\\

Combining these two results with the previous section, we can conclude to the following compactness result.
\begin{cor}\label{cor:conv}{\bf (Compactness)}\\
Assume we  have sequences of functions
  $u_k$ and $P_k$ such that
$$\left\{\begin{array}{l}
\left.\begin{array}{l}
Hu_k=f_k\chi_{\{u_k>0\}}\\
u_k\geq 0
\end{array}\right|\quad \mbox{in}\quad Q_1^-,\\
\ f_k(0)\ge 0,
\end{array}\right. $$
and
$$\left\{\begin{array}{l}
HP_k=f_k(0)\chi_{\{P_k>0\}}\\
P_k\geq 0
\end{array}\right|\quad \mbox{in}\quad Q_1^- .$$
Assume further that with $w_k=u_k-P_k$ there holds for some $p\in (1,+\infty)$
\begin{equation}\label{eq:bound}
  ||w_k||_{L^p(Q_1^-)}+\omega(f_k,1)\leq C.
\end{equation}
Then there is a subsequence of $w_k$ converging in $L^p(Q_\frac{1}{2}^-)$.
Moreover there is $q>1$ such that 
\begin{equation}\label{eq:r102}
||w_k||_{L^q_t(W^{1,q}_x)(Q_{\frac12}^-)}\le C.
\end{equation}
\end{cor}

\noindent {\bf Proof of Corollary \ref{cor:conv}}\\ 
The proof is divided into two parts.\\
\noindent {\bf Part 1: (Convergence a.e.)}\\
From Lemma \ref{lem:Linfty} and
Corollary \ref{cor:w1q} it follows that (for $q>1$)
\begin{equation}\label{eq:bound2}
 ||w_k||_{L_t^\infty(L_x^p)(Q_\frac{1}{2}^-)}, ||\nabla
w_k||_{L_t^q(L_x^q)(Q_\frac{1}{2}^-)}\leq C.
 \end{equation}
Moreover, 
$$
\dd_t w_k=\div(\nabla w_k)+f_k(0)(\chi_{\{u_k>0\}}-\chi_{\{P_k>0\}})+(f_k-f_k(0))\chi_{\{u_k>0\}}=a_k+b_k+c_k.
$$
By Lemma \ref{lem:inclusion}, up to reducing $q<\frac{n}{n-1}$, we have $\nabla L^q(B_r)\subset W^{-1,q}(B_r)$ for any
$B_r$. Then Corollary \ref{cor:w1q}, Lemma \ref{lem:L1} and
\eqref{eq:bound} imply that
$$
||a_k||_{L^1_t(W^{-1,q}_x)(Q_\frac{1}{2}^-)},||b_k||_{L^1(Q_\frac{1}{2}^-)},
||c_k||_{L^1(Q_\frac{1}{2}^-)}\leq C.
$$
Lemma \ref{lem:inclusion} also implies that $L^1(B_r)\subset
W^{-1,q}(B_r)$ for any $B_r$.
Therefore
$$
||\dd_t w_k||_{L^1_t(W^{-1,q}_x)(Q_\frac{1}{2}^-)}\leq C'.
$$
Now we wish to apply Theorem \ref{th:aubin}. Let
$X_0=W^{1,q}(B_\frac12)$, $X=L^1(B_\frac12)$ and $X_1=W^{-1,q}(B_\frac12)$. We can then conclude that there is a subsequence of
$w_k$ that converges in $L_t^1(L_x^1)(Q_\frac{1}{2}^-)=L^1(Q_\frac{1}{2}^-)$. Hence, there is a subsequence that converges
a.e..\\
\noindent {\bf Part 2: ($L^p$-convergence)}\\ 
We wish to apply Theorem
\ref{th:aubin} to the sequence $|w_k|^p$. In order to be able to do
that, we need estimates.\\
\noindent {\bf Step A: Bound on $\nabla |w_k|^p$}\\ 
We set $W_k = w_k|w_k|^{\frac{p}{2}-1}$. We observe that from Lemma \ref{lem:Linfty} and \eqref{eq:bound}, we have
\begin{equation}\label{eq:r104}
||W_k||_{L^\infty_t(L^2_x)(Q_\frac{1}{2}^-)}\leq C_1.
\end{equation}
On the other hand, Lemma \ref{lem:cac1} implies
\begin{equation}\label{eq:r106}
||\nabla W_k||_{L^2_t(L^2_x)(Q_\frac{1}{2}^-)}\leq C_2.
\end{equation}
Therefore, from the elliptic Sobolev embedding (with an abuse of notation if $n=1$)
\begin{equation}\label{eq:r105}
||W_k||_{L^2_t(L^{\frac{2n}{n-2}}_x)(Q_\frac{1}{2}^-)}\leq C_3.
\end{equation}
Then the interpolation between (\ref{eq:r104}) and (\ref{eq:r105}) gives for any $\alpha\in (0,1)$
(this is a classical result which can for instance be easily deduced from the $L^p$-interpolation in Brezis \cite{Brezis} page 57)
$$||W_k||_{L^{p^\alpha_t}_t(L^{p^\alpha_x}_x)(Q_\frac{1}{2}^-)}\leq C_4,$$
with
$$\frac{1}{p^\alpha_t} =   \frac{1-\alpha}{\infty}+ \frac{\alpha}{2}< \frac12,\quad \frac{1}{p^\alpha_x} = \frac{1-\alpha}{2}+ \frac{\alpha}{\left(\frac{2n}{n-2}\right)}< \frac12.$$
This implies the existence of some $p_0>2$ such that
\begin{equation}\label{eq:r107}
||W_k||_{L^{p_0}(Q_\frac{1}{2}^-)}\leq C_5.
\end{equation}
Then, with the $q$ given in Part 1
\begin{align*}
\frac{1}{2^q}\int_{Q_\frac{1}{2}^-}|\nabla |w_k|^p|^q 
& = \frac{1}{2^q}\int_{Q_\frac{1}{2}^-}|2W_k \nabla W_k|^q\\
& \le  \left(\int_{Q_\frac{1}{2}^-}|\nabla W_k|^2\right)^{\frac{q}{2}} 
\left(\int_{Q_\frac{1}{2}^-}|W_k|^{q(\frac{2}{q})'}\right)^{\frac{1}{(\frac{2}{q})'}}.
\end{align*}
But $q(\frac{2}{q})'=\frac{2q}{2-q}$ is increasing in $q$ with value $2$ for $q=1$. Therefore, under our assumptions, and up to reducing $q>1$, we can chose $q$ such that $q(\frac{2}{q})'\le p_0$. We then use (\ref{eq:r106}) and (\ref{eq:r107}) to conclude that
\begin{equation}\label{eq:bound3} 
||\nabla |w_k|^p||_{L^q(Q_\frac{1}{2}^-)}\leq C_6.
\end{equation}
\noindent {\bf Step B: Using the PDE to conclude}\\ 
In order to obtain some information about $\dd_t |w_k|^p$ we need to
play with the equation for $w_k$ again. Multiplication by $|w_k|^{p-1}\sgn
(w_k)$ gives
\begin{align*}
& -\partial_t\left(\frac{|w_k|^p}{p}\right)+|w_k|^{p-1}\sgn (w_k)\lap w  \\
&=\left\{f_k(0)(\chi_{\{u_k>0\}}-\chi_{\{P_k>0\}})+(f_k-f_k(0))\chi_{\{u_k>0\}} \right\}|w_k|^{p-1}\sgn (w_k).
\end{align*}
Rearranging a bit this yields
$$\begin{array}{lll}
-\partial_t\left(\frac{|w_k|^p}{p}\right)&= &\quad \left\{f_k(0)(\chi_{\{u_k>0\}}-\chi_{\{P_k>0\}})+(f_k-f_k(0))\chi_{\{u_k>0\}} \right\}|w_k|^{p-1}\sgn (w_k)\\
\\
&&\quad-\div((|w_k|^{p-1}\sgn w_k\nabla w_k)\\
\\
&&\quad +(p-1)|w_k|^{p-2}|\nabla w_k|^2\\
\\
&=&a_k+b_k+c_k.
\end{array}$$
An observant reader might see that a priori, the calculations above are
not valid other than in some formal sense. However, a simple approximation
argument can make this rigorously justified. Now, Lemma \ref{lem:cac1}, the bound (\ref{eq:bound}) and H\"older's inequality imply that 
\begin{equation}\label{eq:acest}
||a_k||_{L^1(Q_\frac{1}{2}^-)},||c_k||_{L^1(Q_\frac{1}{2}^-)}\leq C.
\end{equation}
Estimate \eqref{eq:bound3} and Lemma \ref{lem:inclusion} imply
$$||b_k||_{L^1_t(W^{-1,q}_x)(Q_\frac{1}{2}^-)}$$
and estimate (\ref{eq:acest}) and Lemma \ref{lem:inclusion} imply
$$
||a_k||_{L^1_t(W^{-1,q}_x)(Q_\frac{1}{2}^-)},||c_k||_{L^1_t(W^{-1,q}_x)(Q_\frac{1}{2}^-)}\leq C.
$$
Hence,
$$
||\dd_t |w_k|^p||_{L^1_t(W^{-1,q}_x)(Q_\frac{1}{2}^-)}\leq C.
$$
Applying Theorem \ref{th:aubin} for the sequence $|w_k|^p$ with $X_0=W^{1,q}(B_\frac12)$, $X=L^1(B_\frac12)$ and $X_1=W^{-1,q}(B_\frac12)$, we find that, up to a subsequence,
$|w_k|^{p}$ converges in $L^1$. This, together with the a.e.-convergence
of a subsequence of $w_k$, implies the existence of a subsequence of
$w_k$ converging in $L^p(Q_\frac{1}{2}^-)$.\\
\noindent {\bf Part 3: (Proof of (\ref{eq:r102}))}\\ 
Finally, (\ref{eq:r102}) follows from Corollary \ref{cor:w1q} and the bound (\ref{eq:bound}).\\
This ends the proof of the corollary.

\section{Decay estimates and the proof of Theorem \ref{th:main}}\label{sec:decay}

The aim of this section is to prove Theorem \ref{th:main}. The key result is Proposition \ref{pro:decay}. We define
$$
N(u,\rho)=\inf_{P\in P_{\textup{reg}}}\left(\frac{1}{\rho^{n+2+2p}}\int_{Q_\rho^-}|u-P|^p\d x\d t\right)^\frac{1}{p},
$$
and
$$
M(u,r)=\sup_{0<\rho\leq r}N(u,\rho),
$$
which is nothing else than the quantity $M_{\textup{reg}}(u,r)$ defined in the introduction.

We will need the following result:

\begin{lem}\label{lem:Mest}{\bf (Estimates of $N$ in larger balls)}\\
If
$$
N(u,1)=||u-P_1||_{L^p(Q_1^-)}, 
$$
with $P_1\in \P_\textup{reg}$, 
then for any $\rho\geq 1$, 
$$
\left(\frac{1}{\rho^{n+2+2p}}\int_{Q_\rho^-}|u-P_1|^p\d x\d t\right)^\frac{1}{p}\leq
  C_1\int_1^{2\rho}\frac{N(u,s)}{s}\d s.
$$
\end{lem}

\noindent {\bf Proof of Lemma \ref{lem:Mest}}\\ 
The proof of this lemma is similar to the proof of Lemma 2.9 in \cite{Mon09}, 
which is proved by decomposing $B_\rho$ into dyadic
balls and estimating the quantities in each of the balls.
We notice in particular that for $\alpha\in [1,2]$, we have
$$N(u,\alpha \rho) \le 2^{2+\frac{n+2}{p}}N(u,2\rho),$$
which is used in order to get the result with the integral of $N$ on the right hand side.\\

\begin{pro}\label{pro:decay}{\bf (Decay estimate)}\\
Let $u$ be a solution of \eqref{eq:obstacle}. Then there are constants $M_0,C_0>0$, $r_0,\lambda,\mu\in (0,1)$
  such that for all $r<r_0$ 
$$
M(u,r)\leq M_0\quad \Longrightarrow \quad M(u,\lambda r)<\mu M(u,r)
\quad\textup{or}\quad 
M(u,r)<C_0\sigma(f,r).
$$
\end{pro}

\noindent {\bf Proof of Proposition \ref{pro:decay}}\\ 
\noindent {\bf Step 1: Construction of sequences and a priori estimates}\\
We argue by contradiction. If this is not true, we can
  find $C_k\to \infty$, $M_k,r_k,\lambda_k\to 0$ and $\mu_k\to 1$ such
  that the statement above fails with the corresponding functions
  $u_k$ and $f_k$, i.e., we have $M(u_k,r_k)\leq M_k$ but
  still 
$$M(u_k,\lambda_k r_k)\geq \mu_k M(u_k,r_k)\quad \mbox{and}\quad M(u_k,r_k)\geq C_k\sigma(f_k,r_k).$$
We note that by our assumption, $M(u_k,\lambda_k r_k)\to 0$. This
implies that we can, passing to another subsequence if possible, assume
that for some $0<\rho_k\leq \lambda_kr_k$
$$
\frac{M(u_k,\lambda_k r_k)}{1+1/k}\leq N(u_k,\rho_k)=\e_k\to 0.
$$
Define the rescaled functions
$$
v_k(x,t)=\frac{u_k(\rho_k x,\rho^2_k t)}{ \rho_k^2}
$$
and
\begin{equation}\label{eq:wk}
w_k(x,t)=\frac{u_k(\rho_k x,\rho^2_k t)-P_k(\rho_kx)}{\e_k \rho_k^2},
\end{equation}
where $P_k\in {\mathcal P}_{\textup{reg}}$ is a half-space function realizing the infimum defining
$N(u_k,\cdot)$ at the level $\rho_k$. Now, we wish to pass to the limit, but before doing that we need to
gather up some informative estimates mainly on the functions $w_k$.

By definition
\begin{equation}\label{eq:N1}
\inf_{P\in \P_\textup{reg}}\left(\int_{Q_1^-}\left|w_k-\frac{P-P_k}{\e_k}\right|^p\d x \d t\right)^\frac{1}{p}=1.
\end{equation}
Moreover, with $s\rho_k\leq r_k$, we have
\begin{align}\label{eq:inf1}
&\displaystyle \inf_{P\in \P_\textup{reg}}\left(\frac{1}{s^{n+2+2p}}\int_{Q_{s}^-}\left|w_k-\frac{P-P_k}{\e_k}\right|^p\d x\d t\right)^\frac{1}{p}
\displaystyle =\frac{N(u_k,s\rho_k)}{\e_k}\nonumber \\
&\displaystyle \leq\frac{ M(u_k,r_k)}{\e_k}\leq  \frac{M(u_k,\lambda_kr_k)}{\e_k\mu_k} \leq \frac{1+1/k}{\mu_k}.
\end{align}
Since
$$
N(v_k,1)=\left(\int_{Q_1^-} |v_k-P_k|^p\d x\d t\right)^\frac{1}{p}, 
$$
we also have for $s\in (1,\frac{r_k}{2\rho_k})$
\begin{equation}\label{eq:subquad}
\begin{array}{ll}
\displaystyle \left(\frac{1}{s^{n+2+2p}}\int_{Q_s^-} |w_k|^p\d x\d
  t\right)^\frac{1}{p}&\displaystyle =\frac{1}{\e_k}\left(\frac{1}{s^{n+2+2p}}\int_{Q_{s}^-} |v_k-P_k|^p\d x\d t\right)^\frac{1}{p}\\
\\
&\displaystyle \leq \frac{C_1}{\e_k}\int_1^{2s}\frac{M(v_k,\tau)}{\tau}\d \tau=
\frac{C_1}{\e_k}\int_1^{2s}\frac{M(u_k,\tau\rho_k)}{\tau}\d \tau\\
\\
&\displaystyle \leq \frac{C_1}{\e_k}\int_1^{2s}\frac{M(u_k,r_k)}{\tau}\d \tau\leq
\frac{C_1}{\e_k} \frac{1+1/k}{\mu_k}N(u_k,\rho_k)\ln 2s\\
\\
&\leq C_2\ln 2s,
\end{array}
\end{equation}
where, in the second line we have used Lemma \ref{lem:Mest}.  Up to rotating our coordinates we can assume that
$$
P_k(x,t)=\frac{1}{2}\left(\max(0,x_1)\right)^2.
$$
Then 
$$H w_k=g_k \textup{ in $\{v_k>0\}\cap \{x_1>0\}$},$$
with 
$$
g_k(x,t)=\frac{f_k(\rho_k x,\rho_k^2 t)-f_k(0)}{\e_k} \quad \mbox{and}\quad f_k(0)=1.
$$
Furthermore, for any $0<r<\frac{r_k}{\rho_k}\to \infty$ as $k\to \infty$, there holds
$$
\left(\frac{1}{|Q_r^-|}\int_{Q_r^-}|g_k|^p\d x\d t\right)^\frac{1}{p}\leq
\frac{\sigma(f_k,\rho_kr)}{\e_k}\leq\frac{M(u_k,r_k)}{C_k\e_k}\leq \frac{M(u_k,\lambda_k
  r_k)}{\mu_kC_k\e_k}\leq \frac{1+1/k}{\mu_kC_k}\to 0.
$$
\noindent {\bf Step 2: Passing to the limit}\\
Corollary \ref{cor:conv} implies that, up to a subsequence, the
sequence $w_k$ converges to some function $w_\infty$ in
$L^p_\loc$. From the $L^p$-bound (\ref{eq:subquad}) on $w_k$ and \eqref{eq:wk} it follows that
$$v_k\to P_\infty=\frac{1}{2}(\max(x_1,0))^2\textup{ in $L^p_\loc$}.
$$
From \eqref{eq:subquad} we deduce that for any $s>1$ there holds
\begin{equation}\label{eq:subquadlimit}
\left(\frac{1}{s^{n+2+2p}}\int_{Q_s^-}|w_\infty|^p\d x\d
t\right)^\frac{1}{p}\leq C_2\ln 2s.
\end{equation}
We define the tangent space of $\P_\textup{reg}$ at $P_\infty$ as
$$
T_{P_\infty}\P_\textup{reg}=\left\{ x_1^+(x,\beta)\textup{ where
    $\beta=(\beta_1,\ldots,\beta_n)\in\R^n$ and $\beta_1=0$} \right\},
$$
which is exactly the set of all possible limits of
$$
\frac{\tilde{P}_k-P_\infty}{\e_k},
$$
as $k\to\infty$, with $\e_k\to 0$ and $\tilde{P}_k\to P_\infty$. Then from \eqref{eq:N1} together with the local $L^p$-convergence
of $w_k$
\begin{equation}\label{eq:distone}
\inf_{q\in T_{P_\infty}\P_\textup{reg}}\left(\int_{Q_1^-}|w_\infty-q|^p\d x\d t\right)^\frac{1}{p}=1.
\end{equation}
Moreover, \eqref{eq:inf1} gives that for any $s>0$, 
\begin{equation}\label{eq:degree2}
\inf_{q\in T_{P_\infty}\P_\textup{reg}}\left(\frac{1}{s^{n+2+2p}}\int_{Q_s^-}|w_\infty-q|^p\d x\d
t\right)^\frac{1}{p}\leq 1.
\end{equation}
From the equation for $w_k$ it follows that 
\begin{equation}\label{eq:wekv}
H w_\infty=0\text{ in } \{P_\infty>0\}=\{x_1>0\}. 
\end{equation}
In addition
$$Hv_k=f_k(\rho_kx,\rho_k^2x)\chi_{\left\{v_k>0\right\}},$$
where the right hand side is bounded in $L^p_{\textup{loc}}$.
Therefore, by the interior estimates and the Sobolev embedding (Theorem
\ref{th:interior} and Theorem \ref{th:sobolev}), for a  subsequence 
$v_k\to P_\infty$ in $L^\infty_\textup{loc}$. Moreover, the weak
non-degeneracy (Corollary \ref{cor:weaknondeg}) implies that for any
compact $K\subset \{x_1<0\}\cap Q_R^-$, we have
$$
v_k\leq C\sigma(f_k,R\rho_k) \text{ in $K$}.
$$
Thus, 
$$
w_k\leq \frac{C\sigma(f_k,R\rho_k)}{\e_k} = C \sigma(g_k,R)\to 0 \quad \text{ in $K$}, 
$$
as $k\to\infty$, and then 
\begin{equation}\label{eq:r120}
w_\infty=0\quad \mbox{in}\quad  \{x_1<0\}.
\end{equation}
Finally, from (\ref{eq:r102}) in Corollary \ref{cor:conv}, we also deduce
$w_\infty\in (L_t^q(W_x^{1,q}))_{\textup{loc}}$ for some $q>1$.\\

\noindent {\bf Step 3: Identification of the limit and contradiction}\\
Since $\{x_1=0\}$ is of codimension one
in space and
does not depend on time, there is a trace of $w_\infty$,
enjoying the estimate
$$
||w_\infty||_{L^q(Q_r^-\cap \{x_1=0\})}\leq
C||w_\infty||_{L_t^q(W_x^{1,q})(Q_r^-\cap \{x_1<0\})},
$$
for any $r>0$. From (\ref{eq:r120}), we deduce that
\begin{equation*}
w_\infty=0\quad \mbox{on}\quad  \{x_1=0\}.
\end{equation*}
Define $w^*$ to be the odd reflection of $w_\infty$ with respect to the
plane $\left\{x_1=0\right\}$. Due to \eqref{eq:wekv}, $w^*$ is caloric in $\R^n\times (-\infty,0)$. Moreover, by
\eqref{eq:subquadlimit}, $w_\infty$  grows strictly slower than
$r^3$ at infinity. This implies that $w^*$ is a caloric polynomial, and
$$w^*(x,t)=\alpha t+P(x),$$ 
where $P$ is a polynomial of degree less than
or equal to two. Taking $s\to 0$ in
\eqref{eq:degree2}, we see that $P$ must be homogeneous of degree two. Indeed, any
linear (spatial) or constant part of $w_\infty$ would make the right
hand side of \eqref{eq:degree2} blow up as $s\to 0$. Since $w^*$
vanishes on $\{x_1=0\}$, we have $\alpha=0$ and $w^*$
must be a spatial harmonic polynomial homogeneous of
degree two. Therefore, as in \cite{Mon09} (Step 3 of the proof of Proposition 6.2), we see that $w^*\in
T_{P_\infty}\P_\textup{reg}$ which contradicts (\ref{eq:distone}).
This ends the proof of the proposition.\\

\noindent {\bf Proof of Theorem \ref{th:main}}\\ 
The theorem follows by combining  Proposition 3.2 in \cite{Mon09} with the present Proposition \ref{pro:decay}.


\section{Acknowledgements}
Both of the authors have partially been supported by the ANR project ''MICA'', grant ANR-08-BLAN-0082. Moreover, the first author has been partially supported by the Chair ''Mathematical modelling and numerical
simulation, F-EADS - Ecole Polytechnique - INRIA'', NTNU and MSRI. 
The authors also thank the MSRI for providing excellent conditions of work and they are very grateful to Peter Lindqvist and Giuseppe Mingione, for the help with finding proper references.


\bibliographystyle{amsplain}
\bibliography{ref.bib}  

\def\cprime{$'$} \def\cprime{$'$}
\providecommand{\bysame}{\leavevmode\hbox to3em{\hrulefill}\thinspace}
\providecommand{\MR}{\relax\ifhmode\unskip\space\fi MR }
\providecommand{\MRhref}[2]{%
  \href{http://www.ams.org/mathscinet-getitem?mr=#1}{#2}
}
\providecommand{\href}[2]{#2}
\begin{thebibliography}{10}

\bibitem{ASU00}
D.~E. Apushkinskaya, H.~Shahgholian, and N.~N. Uraltseva, \emph{Boundary
  estimates for solutions of a parabolic free boundary problem}, Zap. Nauchn.
  Sem. S.-Peterburg. Otdel. Mat. Inst. Steklov. (POMI) \textbf{271} (2000),
  no.~Kraev. Zadachi Mat. Fiz. i Smezh. Vopr. Teor. Funkts. 31, 39--55, 313.

\bibitem{ASU03}
D.~E. Apushkinskaya, N.~N. Ural{\cprime}tseva, and Kh. Shakhgolyan, \emph{On
  the {L}ipschitz property of the free boundary in a parabolic problem with an
  obstacle}, Algebra i Analiz \textbf{15} (2003), no.~3, 78--103.

\bibitem{BDM05}
A.~Blanchet, J.~Dolbeault, and R.~Monneau, \emph{On the one-dimensional
  parabolic obstacle problem with variable coefficients}, Elliptic and
  parabolic problems, Progr. Nonlinear Differential Equations Appl., vol.~63,
  Birkh\"auser, Basel, 2005, pp.~59--66.

\bibitem{Bla06}
Adrien Blanchet, \emph{On the singular set of the parabolic obstacle problem},
  J. Differential Equations \textbf{231} (2006), no.~2, 656--672.

\bibitem{BDM06}
Adrien Blanchet, Jean Dolbeault, and R{\'e}gis Monneau, \emph{On the continuity
  of the time derivative of the solution to the parabolic obstacle problem with
  variable coefficients}, J. Math. Pures Appl. (9) \textbf{85} (2006), no.~3,
  371--414.

\bibitem{Brezis}
Haim Brezis, \emph{Analyse fonctionnelle}, Collection Mathematiques Appliquees
  pour la Maitrise. [Collection of Applied Mathematics for the Master's
  Degree], Masson, Paris, 1983, Theorie et applications. [Theory and
  applications].

\bibitem{CPS04}
Luis Caffarelli, Arshak Petrosyan, and Henrik Shahgholian, \emph{Regularity of
  a free boundary in parabolic potential theory}, J. Amer. Math. Soc.
  \textbf{17} (2004), no.~4, 827--869 (electronic).

\bibitem{EL11}
Anders Edquist and Erik Lindgren, \emph{Regularity of a parabolic free boundary
  problem with h\"older continuous coefficients}, preprint (2011).

\bibitem{evans2009partial}
Lawrence~C. Evans, \emph{Partial differential equations}, Graduate Studies in
  Mathematics, vol.~19, American Mathematical Society, Providence, RI, 1998.

\bibitem{LSU67}
O.~A. Lady{\v{z}}enskaja, V.~A. Solonnikov, and N.~N. Ural{\cprime}ceva,
  \emph{Linear and quasilinear equations of parabolic type}, Translated from
  the Russian by S. Smith. Translations of Mathematical Monographs, Vol. 23,
  American Mathematical Society, Providence, R.I., 1967.

\bibitem{Lie96}
Gary~M. Lieberman, \emph{Second order parabolic differential equations}, World
  Scientific Publishing Co. Inc., River Edge, NJ, 1996.

\bibitem{LM}
Erik Lindgren and R\'egis Monneau, \emph{Pointwise estimates at singular points
  for the parabolic obstacle problem}, work in progress (2011).

\bibitem{Mon09}
R.~Monneau, \emph{Pointwise estimates for {L}aplace equation. {A}pplications to
  the free boundary of the obstacle problem with {D}ini coefficients}, J.
  Fourier Anal. Appl. \textbf{15} (2009), no.~3, 279--335.

\bibitem{Sim87}
Jacques Simon, \emph{Compact sets in the space {$L\sp p(0,T;B)$}}, Ann. Mat.
  Pura Appl. (4) \textbf{146} (1987), 65--96.

\bibitem{W06}
X.J. Wang, \emph{{Schauder Estimates for Elliptic and Parabolic equations}},
  Chinese Annals of Mathematics-Series B \textbf{27} (2006), no.~6, 637--642.

\bibitem{ZC02}
X.~Zou and Y.Z. Chen, \emph{{Fully Nonlinear Parabolic Equations and the Dini
  Condition}}, Acta Mathematica Sinica \textbf{18} (2002), no.~3, 473--480.

\end{thebibliography}

\end{document}